\documentclass[twoside, onecolumn, 10pt]{article}
\usepackage{latexsym}
\usepackage{amssymb}
\usepackage{amsmath}
\usepackage{mathrsfs}
\newtheorem{theorem}{Theorem}[section]
\newtheorem{lemma}[theorem]{Lemma}

\newtheorem{definition}[theorem]{Definition}

\newcommand{\nb}{\nonumber}
\newcommand{\R}{{\mathscr R}}

\begin{document}

\begin{center}
{\Large  Solutions of the matrix inequalities $BXB^{*}\leqslant^{-}\!A$  in the minus partial ordering and $BXB^{*} \leqslant^{{\footnotesize{\rm L}}}\!A$  in the L\"owner partial ordering}
\end{center}
\begin{center}
{\large Yongge Tian}
\end{center}

\begin{center}
{\it {\footnotesize CEMA, Central University of Finance
and Economics, Beijing 100081, China}}
\end{center}

\renewcommand{\thefootnote}{\fnsymbol{footnote}}
\footnotetext{{\it E-mail Address:} yongge.tian@gmail.com}

{\small
\noindent {\bf Abstract.} Two matrices $A$ and $B$ of the same size are said
to satisfy the minus partial ordering, denoted by $B\leqslant^{-}A$,
 iff the rank subtractivity equality ${\rm rank}(\, A - B\,) = {\rm rank}(A) -{\rm rank}(B)$ holds;
 two complex Hermitian matrices $A$ and $B$ of the same size are said
to satisfy the  L\"owner partial ordering, denoted by $B\leqslant^{\rm L} A$,
iff the difference $ A - B$ is nonnegative definite. In this note, we establish general
solution of the inequality $BXB^{*} \leqslant^{-}A$ induced from the minus partial ordering,
and general solution of the inequality $BXB^{*} \leqslant^{\rm L} A$ induced from
the L\"owner partial ordering, respectively, where $(\cdot)^{*}$ denotes the conjugate transpose of a
complex matrix. As consequences, we give closed-form expressions for the shorted matrices of
$A$ relative to the range  of $B$ in the minus  and L\"owner partial orderings, respectively,
and show that these two types of shorted matrices in fact are the same.\\

\noindent {\bf Mathematics Subject Classifications:} 15A03; 15A09; 15A24; 15B57

\medskip

\noindent  {\bf Keywords}: Minus partial ordering;  L\"owner partial ordering; Hermitian matrix;
matrix function; matrix equation; Moore--Penrose inverse; shorted matrix;  rank; inertia
}

\section{Introduction}
\renewcommand{\theequation}{\thesection.\arabic{equation}}
\setcounter{section}{1}
\setcounter{equation}{0}

Throughout this note, let $\mathbb C^{m\times n}$ and  ${\mathbb C}_{{\rm H}}^{m}$ denote the
collections of all $m\times n$ complex matrices and all $m\times m$ complex Hermitian matrices,
respectively; the symbols $A^{*}$, $r(A)$ and ${\mathscr R}(A)$ stand for the conjugate
transpose, the rank and the range (column space) of a matrix $A \in
\mathbb C^{m \times n}$, respectively; $I_m$ denotes the identity
matrix of order $m$; $[\, A, \, B\,]$ denotes a row block matrix
consisting of $A$ and $B$.  The Moore--Penrose inverse of a matrix
$A \in {\mathbb C}^{m \times n}$,
denoted by $A^{\dag}$, is defined to be the unique matrix $X \in
{\mathbb C}^{n \times m}$ satisfying the matrix equations
$$
{\rm (i)} \  AXA = A,  \ \ {\rm (ii)} \ XA X = X, \ \ {\rm (iii)} \
(AX)^{*} = AX, \ \  {\rm (iv)} \ (XA)^{*} =XA.
$$
Further, let $E_A = I_m - AA^{\dag}$ and $F_A = I_n - A^{\dag}A$,
both of which are orthogonal projectors and their ranks are given
by $r(E_A) = m - r(A)$ and $r(F_A) = n - r(A)$. A well-known
property of the Moore--Penrose inverse is $(A^{\dag})^{*} =
(A^{*})^{\dag}$. Hence, if $A = A^{*}$, then both $A^{\dag} =
(A^{\dag})^{*}$ and $AA^{\dag} =A^{\dag}A$ hold.
The inertia of a matrix $A \in {\mathbb C}_{{\rm H}}^{m}$
 is defined to be the triplet ${\rm In}(A)= \{\, i_{+}(A), \, i_{-}(A), \, i_{0}(A)\,\},$
where $i_{+}(A)$, $i_{-}(A)$ and $i_0(A)$ are the numbers of the
positive, negative and zero eigenvalues of $A$ counted with
multiplicities, respectively. For a matrix $A\in
{\mathbb C}_{{\rm H}}^{m}$, both $r(A) = i_{+}(A) + i_{-}(A)$ and
$i_0(A) = m - r(A)$ hold.

The definitions of two well-known partial orderings on  matrices of the same size are given below.

\begin{definition} \label{TS11}
{\rm
\begin{enumerate}
\item[{\rm(a)}]  Two matrices $A, \, B \in \mathbb C^{m \times n}$ are said to satisfy the minus
partial ordering,  denoted by
 $B\leqslant^{-} A$, iff the rank subtractivity equality $r(\, A - B\,) = r(A) - r(B)$ holds,  or
  equivalently, both ${\mathscr R}(A - B) \cap {\mathscr R}(B) = \{0\}$ and
   ${\mathscr R}(A^{*} - B^{*}) \cap {\mathscr R}(B^{*})
   = \{0\}$ hold.

\item[{\rm(b)}] Two matrices $A,  \, B \in {\mathbb C}_{{\rm H}}^{m}$ are said to satisfy the
 L\"owner partial ordering,  denoted by $B\leqslant^{{\footnotesize{\rm L}}} A$,
  iff the difference $A - B$ is nonnegative definite$,$ or equivalently, $A - B = UU^{*}$ for some matrix $U.$
\end{enumerate}
}
\end{definition}

In this note, we consider the following two matrix inequalities
\begin{align}
& BXB^{*} \leqslant^{-} A,
\label{11}
\\
& BXB^{*} \leqslant^{{\footnotesize{\rm L}}} A
\label{12}
\end{align}
 induced from the minus and L\"owner partial orderings, and examine the relations of their solutions,
  where $A \in {\mathbb C}_{{\rm H}}^{m}$ and $B \in \mathbb C^{m \times n}$ are given, and
$X \in {\mathbb C}_{{\rm H}}^{n}$ is unknown. This consideration is motivated by some recent work
on  rank and inertia optimizations of $A -  BXB^{*}$ in \cite{LT,T-laa10,TL}. We shall derive general solutions
of (\ref{11}) and (\ref{12}) by using the given matrices and their generalized inverses,
and then discuss some algebraic properties of these solutions. In particular, we give solutions
of the following constrained rank and L\"owner partial ordering optimization problems
\begin{align}
& \max_{BXB^{*} \leqslant^{-} A} r(BXB^{*}), \ \ \ \ \ \ \ \ \
\min_{BXB^{*} \leqslant^{-} A}  r(\, A -  BXB^{*} \,),
\label{13}
\end{align}\\[-11mm]
\begin{align}
& \max_{\leqslant^{{\footnotesize{\rm L}}}} \{\, BXB^{*} \,| \, BXB^{*} \leqslant^{{\footnotesize{\rm L}}} A \,\},
\ \ \ \ \ \ \ \min_{\leqslant^{{\footnotesize{\rm L}}}} \{\, A - BXB^{*} \,| \, BXB^{*} \leqslant^{{\footnotesize{\rm L}}} A \,\}.
\label{14}
\end{align}
Eqs.\,(\ref{11}) and (\ref{12}) are equivalent to determining elements in the following matrix sets:
 \begin{align}
& {\cal S}_1 = \{\, Z  \in {\mathbb C}_{{\rm H}}^{m} \ | \  Z \leqslant^{-} A, \ \ \R(Z) \subseteq \R(B) \, \},
\label{15}
\\
&  {\cal S}_2 = \{\, Z  \in {\mathbb C}_{{\rm H}}^{m} \ | \ Z \leqslant^{{\footnotesize{\rm L}}} A,
\ \ \R(Z) \subseteq \R(B) \,\}.
\label{16}
\end{align}
The matrices $Z$ in (\ref{15}) and (\ref{16}) can be regarded as two constrained approximations of the matrix
$A$ in partial orderings. In particular, a matrix $Z \in {\cal S}_1$
that has the maximal possible rank is called a shorted matrix of $A$
relative to $\R(B)$ in the minus partial ordering (see \cite{Mi,MP}); while the maximal matrix
in ${\cal S}_2$ is called a shorted matrix of $A$ relative to $\R(B)$ in the L\"owner partial
ordering (see \cite{An,AT}). Our approaches to (\ref{11})--(\ref{14}) link some previous and
recent work in  \cite{An,AT,Go,Gr,HS,Mi,MBS,MP}  on shorted matrices of $A$ relative to given subspaces
in partial orderings, and some recent work on the rank and inertia of the matrix function $A - BXB^{*}$
in \cite{LT,T-laa10,TL}. It is obvious that there always exists a matrix $X$ that satisfies (\ref{11}),
say, $X =0$. Hence,  what we need to do is to derive a general expression
of $X$ that satisfies (\ref{11}). Eq.\,(\ref{12}) may have no solutions unless
the given matrices  $A$ and $B$ in (\ref{12}) satisfy certain conditions.

This note is organized as follows. In Section 2, we present some known results on ranks and inertias
of matrices and matrix equations, and then solve two homogeneous matrix equations with
symmetric patterns. In Section 3, we use the results obtained in Section 2 to derive the
general solution of (\ref{11}), and give an analytical expression for the shorted matrix of
$A$ relative to $\R(B)$ in the minus  partial ordering. In Section 4,
we derive necessary and sufficient conditions for (\ref{12}) to have
a solution,  and then give the general solution of (\ref{12}). We show in Section 5 an interesting
fact that the shorted matrices of $A$ relative to $\R(B)$  in  the minus and L\"owner partial
orderings are the same.

\section{Preliminary results}
\renewcommand{\theequation}{\thesection.\arabic{equation}}
\setcounter{section}{2}
\setcounter{equation}{0}

In order to characterize matrix equalities that involve the Moore--Penrose inverses, we need the
following rank and inertia expansion formulas.

\begin{lemma}[\cite{MS}] \label{T21}
\, Let $A \in \mathbb C^{m \times n}, \, B \in \mathbb C^{m \times
k}$ and $C \in \mathbb C^{l \times n}$ be given$.$  Then$,$ the following
rank expansion formulas hold
\begin{align}
r[\, A, \, B \,] & = r(A) + r(E_AB)  = r(B) + r(E_BA),
\label{21}
\\
r\!\left[\!\begin{array}{c} A \\ C \end{array} \!\right] & = r(A) +
r(CF_A) = r(C) + r(AF_C),
\label{22}
\\
 r\!\left[\!\begin{array}{cccc} A & B
\\
C & CA^{\dag}B
\end{array}\!\right] & =  r\!\left[\!\begin{array}{c} A \\ C \end{array} \!\right] +
r[\, A, \, B \,] - r(A).
\label{qq23}
\end{align}
\end{lemma}

\begin{lemma} [\cite{T-laa10}] \label{T22}
Let $A \in {\mathbb C}_{{\rm H}}^{m},$ $B \in \mathbb C^{m\times n},$
 and $D \in  {\mathbb C}_{{\rm H}}^{n}.$ Then$,$ the following
inertia expansion formulas hold
\begin{align}
&  i_{\pm}\!\left[\!\!\begin{array}{cc}  A  & B  \\ B^{*}  & 0 \end{array}
\!\!\right] = r(B) + i_{\pm}(E_BAE_B),
\label{qq24}
\\
& i_{\pm}\left[\!\!\begin{array}{cc}  A  & B  \\ B^{*}  & D \end{array}
\!\!\right] = i_{\pm}(A) + i_{\pm}(\,  D - B^{*}A^{\dag}B \,)\ \ for \   {\mathscr R}(B) \subseteq {\mathscr R}(A).
\label{qq25}
\end{align}
\end{lemma}

In order to solve  (\ref{11}) and (\ref{12}), we also need the following results on
solvability conditions and general solutions of two simple linear matrix equations.

\begin{lemma}  \label{T23}
Let $A\in \mathbb C^{m\times n}$  and $B\in {\mathbb C}^{m\times p}$ be given$.$ Then$,$ the following hold$.$
\begin{enumerate}
\item[{\rm(a)}] {\rm \cite{Pen}} The matrix equation $AX = B$ is consistent
if and only if ${\mathscr R}(B) \subseteq {\mathscr R}(A).$  In this case$,$ the general
solution can be written as $X = A^{\dag}B + F_AU,$ where $U \in {\mathbb C}^{n\times p}$ is
arbitrary$.$

\item[{\rm(b)}] {\rm \cite{KM}}  Under $B\in \mathbb C^{m\times n},$ the matrix equation $AX = B$
has a solution  $0\leqslant^{{\footnotesize{\rm L}}} X \in {\mathbb C}_{{\rm H}}^{n}$ if and
only if ${\mathscr R}(B) \subseteq {\mathscr R}(A),$ $AB^{*} \geqslant^{{\footnotesize{\rm L}}} 0$
and $r(AB^{*}) =r(B).$ In this case$,$ the
general nonnegative definite solution can be written as
\begin{equation}
X = B^{*}(AB^{*})^{\dag}B   + F_AUF_A,
 \label{26}
\end{equation}
where $0 \leqslant^{{\footnotesize{\rm L}}}  U\in {\mathbb C}_{{\rm H}}^{n}$ is arbitrary$.$
\end{enumerate}
\end{lemma}

\begin{lemma} \label{T24}
Let $A \in \mathbb C^{m \times n}$ and $B \in \mathbb C_{{\rm
H}}^{m}$ be given$.$ Then$,$  the following hold$.$
\begin{enumerate}
\item[{\rm(a)}] {\rm \cite{Gr}} The matrix equation
\begin{equation}
AXA^{*} = B
\label{27}
\end{equation}
has a solution $X \in \mathbb C_{{\rm H}}^{n}$ if and only if
${\mathscr R}(B) \subseteq {\mathscr R}(A),$ or equivalently$,$
$AA^{\dag}B = B.$

\item[{\rm(b)}] {\rm \cite{T-laa10}} Under $X \in \mathbb C_{{\rm H}}^{n},$ the general Hermitian solution of
 {\rm (\ref{27})} can be written in the following two forms
\begin{align}
X & = A^{\dag}B(A^{\dag})^{*} +  U- A^{\dag}A UA^{\dag}A,
\label{28}
\\
X & = A^{\dag}B(A^{\dag})^{*} +  F_AV + V^{*}F_A,
\label{29}
\end{align}
respectively$,$ where $U \in \mathbb C_{{\rm H}}^{n}$ and $V \in {\mathbb C}^{n \times n}$ are arbitrary$.$
\end{enumerate}
\end{lemma}

\begin{lemma}\label{T25}
Let $P \in \mathbb C^{m \times n}$ and $Q \in \mathbb C^{m \times k}$ be given$.$
Then$,$ the general solutions  $X \in {\mathbb C}_{{\rm H}}^{n}$ and
$Y \in {\mathbb C}_{{\rm H}}^{k}$ of the  matrix equation
\begin{align}
PXP^{*} = QYQ^{*}
 \label{210}
\end{align}
can be written as
\begin{align}
X = X_1WX^{*}_1 + X_2,  \ \ \ \ Y = Y_1WY^{*}_1 + Y_2,
 \label{211}
\end{align}
where $W  \in {\mathbb C}_{{\rm H}}^{m}$ is arbitrary$,$ and
$X_1 \in \mathbb C^{n \times m},$ $X_2 \in {\mathbb C}_{{\rm H}}^{n},$ $Y_1  \in \mathbb C^{k \times m}$
and $Y_2 \in {\mathbb C}_{{\rm H}}^{k}$ are the general solutions of the following
 matrix equations
\begin{align}
PX_1 = QY_1, \ \ \ PX_2P^{*} = 0, \ \ \  QY_2Q^{*} = 0,
\label{212}
\end{align}
or alternatively$,$ the general solution of {\rm (\ref{210})} can be written in the
following pair of parametric form
\begin{align}
X & = \widehat{I}_nF_HUF_{H}\widehat{I}^{*}_n
 + U_1 - P^{\dag}PU_1P^{\dag}P,
\label{213}
\\
Y & = \widetilde{I}_k F_HUF_{H}\widetilde{I}^{*}_k  + U_2 - Q^{\dag}QU_2Q^{\dag}Q,
\label{214}
\end{align}
where $H = [\, P, \, -Q \,],$  $\widehat{I}_n = [\, I_n, \, 0\,],$  $\widetilde{I}_k = [\, 0, \, I_k\,],$
and $U \in {\mathbb C}_{{\rm H}}^{n+k},$ $U_1 \in {\mathbb C}_{{\rm H}}^{n}$ and $U_2 \in {\mathbb C}_{{\rm H}}^{k}$
are arbitrary$.$
\end{lemma}

\noindent  {\bf Proof} \, It is easy to verify that the pair of matrices $X$ and $Y$ in (\ref{211})
 are both Hermitian. Substituting the pair of matrices into (\ref{210}) gives
$$
PXP^{*} = PX_1WX^{*}_1P^{*} = QY_1WY^{*}_1Q^{*} = QYQ^{*},
$$
which shows that (\ref{211}) satisfies (\ref{210}). Also, assume that
$X_0$ and $Y_0$ are any pair of solutions of (\ref{210}), and set
$$
W = (PX_0P^{*})^{\dag} = (QY_0Q^{*})^{\dag}, \ \ X_1 = P^{\dag}PX_0P^{*}, \ \
Y_1 = Q^{\dag}QY_0Q^{*},
$$
$$ X_2 = X_0 -  P^{\dag}PX_0P^{\dag}P, \ \ \ \  Y_2 = Y_0 -  Q^{\dag}QY_0Q^{\dag}Q.
$$
Then$,$ (\ref{211}) reduces to
\begin{align*}
X & = P^{\dag}PX_0P^{*}(PX_0P^{*})^{\dag}(P^{\dag}PX_0P^{*})^{*} + X_0 - P^{\dag}PX_0P^{\dag}P
\\
& =  P^{\dag}(PX_0P^{*})(PX_0P^{*})^{\dag}(PX_0P^{*})(P^{\dag})^{*} + X_0 - P^{\dag}PX_0P^{\dag}P
\\
& =  P^{\dag}PX_0P^{\dag}P + X_0 - P^{\dag}PX_0P^{\dag}P = X_0,
\\
Y&  = Q^{\dag}QY_0Q^{*}(QY_0Q^{*})^{\dag}(Q^{\dag}QY_0Q^{*})^{*} + Y_0 - Q^{\dag}QY_0Q^{\dag}Q
\\
& =  Q^{\dag}(QY_0Q^{*})(QY_0Q^{*})^{\dag}(QY_0Q^{*})(Q^{\dag})^{*}+ Y_0 - Q^{\dag}QY_0Q^{\dag}Q
\\
& =  Q^{\dag}QY_0Q^{\dag}Q + Y_0 - Q^{\dag}QY_0Q^{\dag}Q = Y_0,
\end{align*}
that is, any pair of solutions of (\ref{210}) can be represented by
(\ref{211}). Thus, (\ref{211}) is the general solution of (\ref{210}).

Solving the latter two equations in (\ref{212}) by Lemma  \ref{T24}(b) yields the following general solutions
\begin{align}
X_2 = U_1 - P^{\dag}PU_1P^{\dag}P,  \ \ \ \ \ Y_2 = U_2 -
Q^{\dag}QU_2Q^{\dag}Q,
\label{215}
\end{align}
where $U_1 \in {\mathbb C}_{{\rm H}}^{n}$ and $U_2\in {\mathbb C}_{{\rm H}}^{k}$ are arbitrary. To solve the first
equation in (\ref{212}),  we rewrite it as $[\,P, \, -Q\,]\!\left[\! \begin{array}{c} X_1  \\ Y_1
\end{array} \!\right]  =0$. Solving this equation by Lemma \ref{T23}(a) gives the general solution
$\left[\! \begin{array}{c} X_1  \\ Y_1 \end{array} \!\right] =
F_HV_1,$ where $V_1$ is an arbitrary matrix. Hence, the general
expressions of $X_1$ and $Y_1$ can be written as
\begin{align}
X_1 = \widehat{I}_n F_HV_1, \ \ Y_1 = \widetilde{I}_k F_HV_1.
\label{216}
\end{align}
Substituting (\ref{215}) and  (\ref{216}) into (\ref{211}) gives
(\ref{213}) and (\ref{214}). \qquad $\Box$

\begin{lemma}\label{T26}
 Let $B \in \mathbb C^{m \times n}$ and $A \in {\mathbb C}_{{\rm H}}^{m}$ be given$.$
  Then$,$ the general solution $X \in {\mathbb C}_{{\rm H}}^{n}$ of the quadratic matrix equation
\begin{align}
(BXB^{*})A(BXB^{*}) = BXB^{*}
\label{217}
\end{align}
can be expressed in the following parametric form
\begin{align}
X = U(U^{*}B^{*}ABU)^{\dag}U^{*} +  V - B^{\dag}BVB^{\dag}B,
\label{218}
\end{align}
where $U \in \mathbb C^{n \times n}$ and $V \in {\mathbb C}_{{\rm H}}^{n}$ are
arbitrary$.$
\end{lemma}

\noindent {\bf Proof} \, Substituting (\ref{218}) into $BXB^{*}$
gives $BXB^{*} = BU(U^{*}B^{*}ABU)^{\dag}U^{*}B^{*}$. It is easy to
verify by the definition of the Moore--Penrose inverse that
\begin{align*}
(BXB^{*})A(BXB^{*}) & = BU(U^{*}B^{*}ABU)^{\dag}U^{*}B^{*}ABU(U^{*}B^{*}ABU)^{\dag}U^{*}B^{*}
\\
& =BU(U^{*}B^{*}ABU)^{\dag}U^{*}B^{*}  = BXB^{*}.
\end{align*}
Hence, (\ref{218}) satisfies (\ref{217}). On the other hand, for
any Hermitian solution $X_0$ of (\ref{217}), set $U
=B^{\dag}BX_0B^{\dag}B$ and $V = X_0$ in (\ref{218}). Then$,$
(\ref{218}) reduces to
\begin{align*}
X & = B^{\dag}BX_0B^{\dag}B(B^{\dag}BX_0B^{*}ABX_0B^{\dag}B)^{\dag}B^{\dag}BX_0B^{\dag}B + X_0 - B^{\dag}BX_0B^{\dag}B
\\
& = B^{\dag}BX_0B^{\dag}B(B^{\dag}BX_0B^{\dag}B)^{\dag}B^{\dag}BX_0B^{\dag}B + X_0 - B^{\dag}BX_0B^{\dag}B
\\
& = B^{\dag}BX_0B^{\dag}B + X_0 - B^{\dag}BX_0B^{\dag}B = X_0.
\end{align*}
This result indicates that all solutions of (\ref{217}) can be
represented through (\ref{218}). Hence, (\ref{218}) is the general
solution of (\ref{217}). \qquad $\Box$

\section{General solution of $BXB^{*} \leqslant^{-} A$}
\renewcommand{\theequation}{\thesection.\arabic{equation}}
\setcounter{section}{3}
\setcounter{equation}{0}

A well-known necessary and sufficient condition for the rank subtractivity equality in Definition 1.1
to hold is
\begin{equation}
r(\, A - B\,) = r(A) - r(B) \Leftrightarrow  \R(B) \subseteq \R(A), \
 \R(B^{*}) \subseteq \R(A^{*})  \ \text{and}  \ BA^{\dag}B = B,
\label{31}
\end{equation}
see \cite{MS}. Applying (\ref{31}) to  (\ref{11}), we can convert (\ref{11}) to a system of matrix equations.

\begin{lemma}\label{T31}
Eq.~{\rm (\ref{11})} is equivalent to the following system of
matrix equations
\begin{align}
 BXB^{*} = AYA,  \ \  (BXB^{*})A^{\dag}(BXB^{*}) = BXB^{*},
 \label{32}
\end{align}
where $Y\in {\mathbb C}_{{\rm H}}^{m}$ is an unknown matrix$.$
\end{lemma}

\noindent {\bf Proof} \, From (\ref{31}), the minus partial order $BXB^{*} \leqslant^{-} A$ in (\ref{11}) is
equivalent to
\begin{align}
\R(BXB^{*})  \subseteq  \R(A) \ \  {\rm and} \ \   (BXB^{*})A^{\dag}(BXB^{*}) = BXB^{*}.
 \label{33}
\end{align}
By Lemma \ref{T24}(a), the first range inclusion in (\ref{33}) holds
if and only if the first matrix equation in  (\ref{32}) is solvable
for $Y$. Thus,  (\ref{32}) and  (\ref{33}) are
equivalent. \qquad $\Box$

\begin{theorem}\label{T32}
Let $A \in {\mathbb C}_{{\rm H}}^{m}$ and $B \in \mathbb C^{m\times n}$ be given$,$
and ${\cal S}_1$ be as given in {\rm (\ref{16})}$.$ Also define
$$
M = \left[\! \begin{array}{cc}  A  & B  \\ B^{*} & 0 \end{array}\!\right]\!, \ \
H = [\, B, \, -A \,], \ \
 \widehat{I}_n = [\, I_n, \, 0\,], \ \ \widehat{B} =[\, B, \, 0\,], \ \ A_1 = E_BA, \ \ B_1 = E_AB.
$$
Then$,$  the following hold$.$
\begin{enumerate}
\item[{\rm(a)}] The general Hermitian solution of the inequality
\begin{equation}
BXB^{*}  \leqslant^{-} A
\label{34}
\end{equation}
can be written as
\begin{align}
& X = \widehat{I}_nF_HU(U^{*}F_{H}\widehat{B}^{*}A^{\dag}\widehat{B}F_HU)^{\dag}U^{*}F_{H}\widehat{I}^{*}_n
+ V - B^{\dag}BVB^{\dag}B,
\label{35}
\end{align}
where  $U\in {\mathbb C}^{(m +n)\times (m +n)}$ and $V \in {\mathbb C}_{{\rm H}}^{n}$  are arbitrary$.$

\item[{\rm(b)}]  The general expression of the matrices in {\rm (\ref{15})} can be written as
\begin{align}
 Z = \widehat{B}F_HU(U^{*}F_{H}\widehat{B}^{*}A^{\dag}\widehat{B}F_HU)^{\dag}
U^{*}F_{H}\widehat{B}^{*},
\label{36}
\end{align}
where $U\in {\mathbb C}^{(m +n)\times (m +n)}$ is arbitrary$.$ The global maximal and minimal inertias and ranks
of $Z$ in {\rm (\ref{36})} and the corresponding $A -Z$  are given by
\begin{align}
& \max_{Z\in {\cal S}_1} i_{\pm}(Z)  =  i_{\mp}(M) + i_{\pm}(A) - r[\, A,\,  B \,],
\label{37}
\\
& \max_{Z\in {\cal S}_1} r(Z)  = r(M) + r(A) - 2r[\, A,\,  B \,],
\label{38}
\\
& \min_{Z\in {\cal S}_1}  i_{\pm}( \, A - Z \, ) = r[\, A,\,  B \,] -  i_{\mp}(M),
\label{39}
\\
& \min_{Z\in {\cal S}_1} r( \, A - Z \, ) = 2r[\, A,\,  B \,] - r(M).
\label{310}
\end{align}
The shorted matrix of $A$ relative to $\R(B),$ denoted by $\phi^{-}(\, A \,|\, B\,),$
which is a matrix $Z$ that satisfies {\rm (\ref{38})}$,$ is given by
\begin{align}
\phi^{-}(\, A \,|\, B\,) = \widehat{B}F_H(F_{H}\widehat{B}^{*}A^{\dag}\widehat{B}F_H)^{\dag}F_{H}\widehat{B}^{*}.
\label{311}
\end{align}
\end{enumerate}
\end{theorem}

\noindent {\bf Proof} \, Applying Lemma \ref{T25} to the first equation in (\ref{32}),
we obtain the general solutions of $X$ and $Y$ as follows
\begin{align}
X = \widehat{I}_nF_HTF_{H}\widehat{I}^{*}_n +  V - B^{\dag}BVB^{\dag}B, \ \
Y = \widehat{I}_mF_HTF_{H}\widehat{I}^{*}_m +  W - A^{\dag}AWA^{\dag}A,
\label{312}
\end{align}
where $T \in {\mathbb C}_{{\rm H}}^{m+n}$, $V  \in {\mathbb C}_{{\rm H}}^{n}$ and
$W \in {\mathbb C}_{{\rm H}}^{m}$ are arbitrary$.$ Substituting (\ref{312}) into the second equation
in (\ref{32}) leads to the following quadratic matrix equation
$$
(\widehat{B}F_HTF_{H}\widehat{B}^{*})A^{\dag}(\widehat{B}F_HTF_{H}\widehat{B}^{*})
= \widehat{B}F_HTF_{H}\widehat{B}^{*}.
$$
By Lemma \ref{T26}, the general solution of this quadratic matrix equation is given by
\begin{align*}
T& = U(U^{*}F_{H}\widehat{B}^{*}A^{\dag}\widehat{B}F_HU)^{\dag}U^{*} +
 W_1 - (\widehat{B}F_H)^{\dag}(\widehat{B}F_H)W_1(\widehat{B}F_H)^{\dag}(\widehat{B}F_H),
\end{align*}
where $U\in \mathbb C^{(m+n) \times (m+n)}$ and $W_1 \in {\mathbb C}_{{\rm H}}^{m+n}$ are arbitrary.
Substituting this $T$ into  the matrix $X$ in (\ref{312}) gives
\begin{align}
X & =  \widehat{I}_nF_HU(U^{*}F_{H}\widehat{B}^{*}A^{\dag}\widehat{B}F_HU)^{\dag}U^{*}F_{H}\widehat{I}^{*}_n
\nonumber
\\
&  \ \ + [\, \widehat{I}_nF_HW_1F_{H}\widehat{I}^{*}_n -
\widehat{I}_nF_H(\widehat{B}F_H)^{\dag}(\widehat{B}F_H)W_1(\widehat{B}F_H)^{\dag}(\widehat{B}F_H)F_{H}\widehat{I}^{*}_n \,]
  + V - B^{\dag}BVB^{\dag}B.
\label{313}
\end{align}
It is easy to verify from $B\widehat{I}_n = \widehat{B}$ that
\begin{align*}
& B[\, \widehat{I}_nF_HW_1F_{H}\widehat{I}^{*}_n -
\widehat{I}_nF_H(\widehat{B}F_H)^{\dag}(\, \widehat{B}F_H)W_1(\widehat{B}F_H \,)^{\dag}
(\widehat{B}F_H)F_{H}\widehat{I}^{*}_n\,]B^{*}
\\
& =\widehat{B}F_HW_1F_{H}\widehat{B}^{*} -
(\widehat{B}F_H)(\widehat{B}F_H)^{\dag}(\, \widehat{B}F_H)W_1(\widehat{B}F_H \,)^{\dag}
(\widehat{B}F_H)(\widehat{B}F_{H})^{*} =0.
\end{align*}
This fact shows that the second term on the right-hand side of (\ref{313}) is a solution to $BXB^{*} =0$. Also, note
 from Lemma \ref{T24}(b) that $V - B^{\dag}BVB^{\dag}B$ is the general solution to $BXB^{*} =0$.
 Hence, the second term on the right-hand side of (\ref{313}) can be represented by the third term of the same side,
 so that (\ref{313}) reduces to  (\ref{35}).

Substituting  (\ref{35}) into $BXB^{*}$
gives
\begin{align}
Z = BXB^{*} = \widehat{B}F_HU(\,U^{*}F_{H}\widehat{B}^{*}A^{\dag}\widehat{B}F_HU\,)^{\dag}
U^{*}F_{H}\widehat{B}^{*},
\label{314}
\end{align}
as required for (\ref{36}). Note further that this $Z$ satisfies
\begin{align}
(U^{*}F_{H}\widehat{B}^{*}A^{\dag})Z(A^{\dag}\widehat{B}F_HU)  = U^{*}F_{H}\widehat{B}^{*}A^{\dag}\widehat{B}F_HU.
\label{315}
\end{align}
Both  (\ref{314}) and  (\ref{315}) imply
\begin{align}
i_{\pm}(Z)= i_{\pm}(U^{*}F_{H}\widehat{B}^{*}A^{\dag}\widehat{B}F_HU) \leqslant
i_{\pm}(\,F_{H}\widehat{B}^{*}A^{\dag}\widehat{B}F_H\,)
\label{316}
\end{align}
and
\begin{align}
\max_{Z\in {\cal S}_1}i_{\pm}(Z) = i_{\pm}(F_{H}\widehat{B}^{*}A^{\dag}\widehat{B}F_H), \ \ \
\max_{Z\in {\cal S}_1}r(Z) = r(F_{H}\widehat{B}^{*}A^{\dag}\widehat{B}F_H).
\label{317}
\end{align}
Recall that the inertia of a Hermitian matrix does not change under Hermitian congruence operations.
 Applying (\ref{qq24}) to $F_{H}\widehat{B}^{*}A^{\dag}\widehat{B}F_H$ and simplifying by
  Hermitian congruence operations, we obtain
\begin{align}
i_{\pm}(F_{H}\widehat{B}^{*}A^{\dag}\widehat{B}F_H) & =
i_{\pm}\!\left[\! \begin{array}{cc} \widehat{B}^{*}A^{\dag}\widehat{B} & H^{*}
\\
H &  0 \end{array} \!\right] - r(H) \nonumber
\\
& = i_{\pm}\!\left[\! \begin{array}{cccc} B^{*}A^{\dag}B & 0  & B^{*}
\\
0 & 0 & -A
\\
B & -A & 0  \end{array} \!\right] - r[\, A, \, B \,]  \nonumber
\\
& = i_{\pm}\!\left[\! \begin{array}{cccc} 0 & \frac{1}{2}B^{*}A^{\dag}A  & B^{*}
\\
\frac{1}{2}A^{\dag}A B  & 0 & -A
\\
B & -A & 0  \end{array} \!\right] - r[\, A, \, B \,] \nb
\\
&  = i_{\pm}\!\left[\! \begin{array}{cccc} 0 & 0  & B^{*}
\\
0 & A & -A
\\
B & -A & 0  \end{array} \!\right] - r[\, A, \, B \,] \nonumber
\\
& = i_{\pm}\!\left[\! \begin{array}{cccc} 0 & 0  & B^{*}
\\
0 & A & 0
\\
B & 0 & -A  \end{array} \!\right] - r[\, A, \, B \,] \nonumber
\\
& = i_{\mp}\!\left[\! \begin{array}{cccc} A  & B
\\
B^{*} & 0  \end{array} \!\right] + i_{\pm}(A)  - r[\, A, \, B \,].
\label{318}
\end{align}
Substituting  (\ref{318}) into  (\ref{317}) leads to  (\ref{37}) and (\ref{38}).  Also, note that
$$
\min_{X \in {\cal S}_1} i_{\pm}(\, A - Z\,)
= i_{\pm}(A) - \max_{X \in {\cal S}_1} i_{\pm}(Z).
$$
Thus, (\ref{39}) and  (\ref{310}) follow from  (\ref{37}) and (\ref{38}).  \qquad $\Box$

\section{General solution of $BXB^{*} \leqslant^{{\footnotesize{\rm L}}} A$}
\renewcommand{\theequation}{\thesection.\arabic{equation}}
\setcounter{section}{4}
\setcounter{equation}{0}

In this section, we derive an analytical expression for the general solution of (\ref{12}) by
using generalized inverses of matrices, and show some algebraic properties of the solution.

\begin{theorem} \label{T41}
Let $A \in {\mathbb C}_{{\rm H}}^{m}$ and $B \in {\mathbb C}^{m\times n}$ be given$,$
  and let ${\cal S}_2$ be as given in {\rm (\ref{16})}$.$
Then$,$  the following hold$.$
\begin{enumerate}
\item[{\rm(a)}]  There exists an $X \in {\mathbb C}^{n}_{{\rm H}}$ such that
\begin{equation}
BXB^{*}  \leqslant^{{\footnotesize{\rm L}}} A
\label{41}
\end{equation}
if and only if
\begin{equation}
E_BAE_B \geqslant^{{\footnotesize{\rm L}}} 0 \ \ \ and  \ \ \  r(E_BAE_B)= r(E_BA),
\label{42}
\end{equation}
or equivalently$,$
\begin{equation}
i_{+}\!\left[\!\! \begin{array}{cc} A & B \\ B^{*} & 0 \end{array}\!\!\right] = r[\, A, \, B\,] \ \ \
and  \ \ \ i_{-}\!\left[\!\! \begin{array}{cc} A & B \\ B^{*} & 0 \end{array}\!\!\right] = r(B).
\label{43}
\end{equation}
 In this case$,$ the general Hermitian solution of {\rm (\ref{41})} can be written in
 the following parametric form
\begin{align}
& X  =  B^{\dag}A(B^{\dag})^{*} - B^{\dag}AE_B(E_BAE_B)^{\dag}E_BA(B^{\dag})^{*} - UU^{*}
 + F_BV + V^{*}F_B,
\label{44}
\end{align}
where $U, \ V \in {\mathbb C}^{n \times n}$ are arbitrary$.$  Correspondingly$,$  the general expression of
the matrices in ${\cal S}_2$ can be written as
\begin{align}
Z  =  A - AE_B(E_BAE_B)^{\dag}E_BA  - BUU^{*}B^{*}.
\label{45}
\end{align}

\item[{\rm(b)}] Under {\rm (\ref{42})}$,$ the shorted matrix of $A$ relative to $\R(B),$ denoted
by $\phi^{{\rm L}}(\, A \,|\, B\,),$ which is the maximizer in ${\cal S}_2,$
can uniquely be written as
\begin{align}
\phi^{{\rm L}}(\, A \,|\, B\,) = A - AE_B(E_BAE_B)^{\dag}E_BA.
\label{46}
\end{align}
The rank and inertia of $\phi^{{\rm L}}(\, A \,|\, B\,)$ and $A - \phi^{{\rm L}}(\, A \,|\, B\,)$ satisfy
\begin{align}
& i_{+}[\,\phi^{{\rm L}}(\, A \,|\, B\,) \,] = i_{+}(A) + r(B) - r[\, A, \, B \,], \ \
\label{47}
\\
& i_{-}[\,\phi^{{\rm L}}(\, A \,|\, B\,) \,]   = i_{-}(A),
\label{8}
\\
& i_{+}[\, A - \phi^{{\rm L}}(\, A \,|\, B\,) \,]  =
 r[\, A - \phi^{{\rm L}}(\, A \,|\, B\,) \,]  = r[\, A, \, B\,] - r(B).
\label{49}
\end{align}
\end{enumerate}
\end{theorem}

\noindent {\bf Proof} \, It is obvious that (\ref{41}) is equivalent to
\begin{equation}
BXB^{*} = A - YY^{*}
\label{410}
\end{equation}
for some matrix $Y$.  In other words, (\ref{41}) can be relaxed to a matrix equation with two
unknown matrices. From Lemma \ref{T24}(a), (\ref{410}) is solvable for $X \in {\mathbb C}^{n}_{{\rm H}}$
 if and only if $E_B(\,A - YY^{*} \,) = 0,$ that is,
\begin{equation}
E_BYY^{*} = E_BA.
\label{411}
\end{equation}
From Lemma \ref{T23}(b), (\ref{411}) is solvable for $YY^{*}$ if and only if
$E_BAE_B \geqslant^{{\footnotesize{\rm L}}} 0$
and $r(E_BAE_B)= r(E_BA)$, establishing (\ref{42}), which is further equivalent to (\ref{43})
by (\ref{21}) and (\ref{qq24}). In this case, the general nonnegative definite solution of
(\ref{411}) can be written as
\begin{equation}
YY^{*} =  AE_B(E_BAE_B)^{\dag}E_BA +  BB^{\dag}WBB^{\dag},
\label{412}
\end{equation}
where $0 \leqslant^{{\footnotesize{\rm L}}} W \in {\mathbb C}_{{\rm H}}^{m}$ is
arbitrary. Substituting the $YY^{*}$ into (\ref{410}) gives
\begin{equation}
BXB^{*} = A -AE_B(E_BAE_B)^{\dag}E_BA -  BB^{\dag}WBB^{\dag}.
\label{413}
\end{equation}
By Lemma \ref{T24}(b), the general Hermitian solution of (\ref{413}) can be written as
\begin{equation}
X = B^{\dag}A(B^{\dag})^{*} - B^{\dag}AE_B(E_BAE_B)^{\dag}E_BA(B^{\dag})^{*} - B^{\dag}W(B^{\dag})^{*}
 + F_BV + V^{*}F_B,
\label{414}
\end{equation}
where $V \in \mathbb C^{n \times n}$ is arbitrary$.$ Replacing the matrix
$0 \leqslant^{{\footnotesize{\rm L}}} B^{\dag}W(B^{\dag})^{*} \in {\mathbb C}_{{\rm H}}^{n}$ in (\ref{414})
 with a general matrix $0 \leqslant^{{\footnotesize{\rm L}}} U \in {\mathbb C}_{{\rm H}}^{n}$ yields
 (\ref{44}), which is also the general Hermitian solution of (\ref{41}). Substituting (\ref{44}) into
 $BXB^{*}$ gives (\ref{45}).

 Eq.\,(\ref{46})  follows from  (\ref{45}) by noticing $BUU^*B^{*} \geqslant^{{\footnotesize{\rm L}}} 0$.

It follows from (\ref{42}) that $\R(E_BAE_B)  = \R(E_BA)$. In this case, applying (\ref{qq25}) to (\ref{46})
and simplifying  by Hermitian congruence transformations, we obtain
\begin{align*}
i_{\pm}[\,\phi^{{\rm L}}(\, A \,|\, B\,)\,] & =
i_{\pm}[\,  A - AE_B(E_BAE_B)^{\dag}E_BA \,]
\\
& = i_{\pm}\!\left[\!\begin{array}{cc} E_BAE_B &  E_BA
\\
AE_B & A \end{array}\!\right] - i_{\pm}(E_BAE_B)
\\
& = i_{\pm}\!\left[\!\begin{array}{cc} 0 & 0 \\  0 & A \end{array}\!\right] - i_{\pm}(E_BAE_B) \nonumber
\\
& = i_{\pm}(A) - i_{\pm}(E_BAE_B),
\\
i_{\pm}[\, A - \phi^{{\rm L}}(\, A \,|\, B\,)\,]  &= i_{\pm}[\,AE_B(E_BAE_B)^{\dag}E_BA\,]
\\
& = i_{\pm}\!\left[\!\begin{array}{cc} - E_BAE_B &  E_BA
\\
AE_B & 0 \end{array}\!\right] - i_{\mp}(E_BAE_B)
\\
&  = i_{\pm}\!\left[\!\begin{array}{cc} 0 &  E_BA
\\
AE_B & 0 \end{array}\!\right] - i_{\mp}(E_BAE_B)
\\
&  = r(E_BA) - i_{\mp}(E_BAE_B).
\end{align*}
Hence, we further find from that (\ref{21}) and (\ref{42}) that
\begin{align*}
 i_{+}[\phi^{{\rm L}}(\, A \,|\, B\,)] &  = i_{+}(A) -  i_{+}(E_BAE_B) = i_{+}(A) -  r(E_BA)  = i_{+}(A) + r(B) - r[\, A, \, B \,],
 \\
i_{-}[\,\phi^{{\rm L}}(\, A \,|\, B\,)\,] &  = i_{-}(A) -i_{-}(E_BAE_B)  = i_{-}(A),
\\
i_{+}[\,A - \phi^{{\rm L}}(\, A \,|\, B\,)\,] &  = r(E_BA) - i_{-}(E_BAE_B) =  r(E_BA)  =
 r[\, A, \, B \,] - r(B),
\\
i_{-}[\,A - \phi^{{\rm L}}(\, A \,|\, B\,)\,] &  = r(E_BA) - i_{+}(E_BAE_B) =0,
\end{align*}
establishing (\ref{47})--(\ref{49}).
\qquad $\Box$

\section{An equality for the shorted matrices of $A$ relative to $\R(B)$ in the minus and L\"owner partial orderings}

\renewcommand{\theequation}{\thesection.\arabic{equation}}
\setcounter{section}{5}
\setcounter{equation}{0}

Since ${\cal S}_1$ and ${\cal S}_2$ in  (\ref{15}) and  (\ref{16}) are defined
from different matrix inequalities, the two sets are not necessarily the same, as demonstrated
in Theorems \ref{T32}(b) and \ref{T41}(a).  However, they may have some common matrices.
In this section, we show an interesting fact that the shorted matrices
of $A$ relative to $\R(B)$ in the minus and L\"owner partial orderings are the same.

\begin{theorem} \label{T51}
Let $A \in {\mathbb C}_{{\rm H}}^{m}$ and $B \in \mathbb C^{m\times n}$ be given$,$  and
${\cal S}_1$ and ${\cal S}_2$  be as given in {\rm (\ref{15})} and  {\rm (\ref{16})}$.$
If {\rm (\ref{41})} has a solution$,$ then the two shorted matrices in  ${\cal S}_1$ and ${\cal S}_2$
are the same, namely$,$
\begin{align}
\phi^{-}(\, A \,|\, B\,)  = \phi^{{\rm L}}(\, A \,|\, B\,).
\label{51}
\end{align}
\end{theorem}

\noindent {\bf Proof} \, Note from (\ref{311}) and (\ref{46}) that (\ref{51}) holds  if and only if
\begin{align}
 \widehat{B}F_H(F_{H}\widehat{B}^{*}A^{\dag}\widehat{B}F_H)^{\dag}F_{H}\widehat{B}^{*} =
  A - AE_B(E_BAE_B)^{\dag}E_BA.
\label{52}
\end{align}
It is easy to derive from (\ref{22}) that
\begin{align}
r(\,\widehat{B}F_H \,) = r\!\left[\!\begin{array}{c} B \\ H \end{array} \!\right] - r(H)
= r(A) + r(B) - r[\, A, \, B \,].
\label{53}
\end{align}
 Under (\ref{42}), (\ref{317}) reduces to
\begin{align}
r(\,F_{H}\widehat{B}^{*}A^{\dag}\widehat{B}F_H \,)  = r(M) + r(A) - 2r[\, A, \, B \,] =  r(A) + r(B) -
 r[\, A, \, B \,].
\label{54}
\end{align}
Both (\ref{53}) and (\ref{54}) imply that
$\R(\,F_{H}\widehat{B}^{*}\,) = \R(\,F_{H}\widehat{B}^{*}A^{\dag}\widehat{B}F_H \,).$  In this case,
applying (\ref{qq25}) to the difference of  both sides of (\ref{52}) and simplifying by elementary
matrix operations, we obtain
\begin{align*}
& r[\,   A - AE_B(E_BAE_B)^{\dag}E_BA  -
\widehat{B}F_H(F_{H}\widehat{B}^{*}A^{\dag}\widehat{B}F_H)^{\dag}F_{H}\widehat{B}^{*} \,]  \nb
\\
& = r\!\left[\!\begin{array}{cc} F_{H}\widehat{B}^{*}A^{\dag}\widehat{B}F_H
& F_{H}\widehat{B}^{*}  \\
\widehat{B}F_H &   A - AE_B(E_BAE_B)^{\dag}E_BA   \end{array}\!\right]
- r(F_{H}\widehat{B}^{*}A^{\dag}\widehat{B}F_H)
\\
& = r\!\left[\!\begin{array}{ccc} \widehat{B}^{*}A^{\dag}\widehat{B}
& \widehat{B}^{*}  & H^{*} \\
\widehat{B} &  A - AE_B(E_BAE_B)^{\dag}E_BA  & 0
\\
H & 0 & 0  \end{array}\!\right] - 2r(H) - r(A) - r(B) + r[\, A, \, B \,] \ \  \mbox{(by (\ref{qq24}))}
\\
& = r\!\left[\!\begin{array}{cccc} B^{*}A^{\dag}B & 0 & B^{*}  & B^{*}
\\
0  & 0  & 0 & -A
\\
B & 0& A - AE_B(E_BAE_B)^{\dag}E_BA  & 0
\\
B & -A & 0 & 0  \end{array}\!\right]  - r(A) - r(B) - r[\, A, \, B \,]
\\
& = r\!\left[\!\begin{array}{cccc} B^{*}A^{\dag}B & 0 & B^{*}  & 0
\\
0  & 0  & 0 & -A
\\
B & 0& A - AE_B(E_BAE_B)^{\dag}E_BA  & 0
\\
0 & -A & 0 & 0  \end{array}\!\right]  - r(A) - r(B) - r[\, A, \, B \,]
\\
& = r\!\left[\!\begin{array}{cccc} B^{*}A^{\dag}B & B^{*}
\\
B & A - AE_B(E_BAE_B)^{\dag}E_BA  \end{array}\!\right] + r(A) - r(B) - r[\, A, \, B \,]
\\
& = r\left(\!\left[\!\begin{array}{cccc} B^{*}A^{\dag}B & B^{*}
\\
B & A \end{array}\!\right] - \!\left[\!\begin{array}{cccc} 0 & 0
\\
0 &  AE_B(E_BAE_B)^{\dag}E_BA  \end{array}\!\right] \right) + r(A) - r(B) - r[\, A, \, B \,]
\\
& = r\!\left[\!\begin{array}{cccc} B^{*}A^{\dag}B & B^{*} & 0
\\
B & A  &  AE_B
\\
0 &  E_BA &  E_BAE_B
\end{array}\!\right]  - r(E_BAE_B)  + r(A) - r(B) - r[\, A, \, B \,] \ \  \mbox{(by (\ref{qq25}))}
\\
& = r\!\left[\!\begin{array}{cccc} B^{*}A^{\dag}B & B^{*} & 0
\\
B & A  & 0
\\
0 &  0 &  0
\end{array}\!\right] - r(E_BA)  +  r(A) - r(B) - r[\, A, \, B \,]
\\
& = r\!\left[\!\begin{array}{cccc} B^{*}A^{\dag}B & B^{*}
\\
B & A
\end{array}\!\right]  +  r(A) - 2r[\, A, \, B \,] = 0 \ \ \ \  \mbox{(by (\ref{qq23}))},
\end{align*}
which means that (\ref{52}) is an equality.  \qquad $\Box$

\medskip

The minus and L\"owner partial orderings in Definition 1.1 can accordingly be defined for linear operators
on a Hilbert space. Also, note that the results in this note are derived from some ordinary algebraic
operations of the given matrices and their Moore--Penrose inverses. Hence, it is no doubt that most of the
conclusions in this note can be extended to operator algebra, in which the Moore--Penrose
inverses of linear operators were defined.


\begin{thebibliography}{}
\def\itemsep{-1ex}
         \small

\bibitem{An} W.N. Anderson, Jr., Shorted operators, SIAM J. Appl. Math. 20(1971), 522--525.

\bibitem{AT} W.N. Anderson, Jr. and G.E. Trapp, Shorted operators II, SIAM J. Appl. Math. 28(1975), 60--71.

\bibitem{Go} H. Goller, Shorted operators and rank decomposition
matrices, Linear Algebra Appl. 81(1986), 207--236.

\bibitem{Gr} J. Gro\ss, A note on the rank-subtractivity ordering,  Linear Algebra
 Appl. 289(1999), 151--160.

\bibitem{HS}  R.E. Hartwig and G.P.H. Styan, On some characterizations of the
"star" partial ordering for matrices and rank subtractivity, Linear Algebra
Appl. 82(1986), 145--16l.

\bibitem{KM}  C.G. Khatri and S.K. Mitra, Hermitian and nonnegative
definite solutions of linear matrix equations, SIAM J. Appl. Math. 31(1976), 579--585.


\bibitem{LT} Y. Liu and Y. Tian, Hermitian-type of singular value decomposition for a pair of matrices
 and its applications, Numer. Linear Algebra Appl.,  DOI: 10.1002/nla.1825.

\bibitem{MS} G. Marsaglia and G.P.H. Styan, Equalities  and  inequalities
for ranks of matrices,  Linear and Multilinear Algebra 2(1974), 269--292.

\bibitem{Mi} S.K. Mitra, The minus partial  order and the  shorted
matrix,  Linear Algebra Appl. 83(1986), 1--27.

\bibitem{MBS} S.K. Mitra, P Bhimasankaram, S.B. Malik,  Matrix Partial Orders,
Shorted Operators and Applications, World Scientific Publishing, 2010.

\bibitem{MP} S.K. Mitra and M.L. Puri, Shorted  matrices--an extended
     concept and some applications, Linear Algebra Appl. 42(1982), 57--79.

\bibitem{Pen} R. Penrose, A generalized inverse for matrices,
Proc. Cambridge Philos. Soc. 51(1955), 406--413.

\bibitem{T-laa10} Y. Tian, Equalities and inequalities for inertias of
Hermitian matrices with applications, Linear Algebra Appl. 433(2010), 263--296.


\bibitem{TL} Y. Tian and Y. Liu, Extremal ranks of some
symmetric matrix expressions with applications, SIAM J. Matrix
Anal. Appl. 28(2006), 890--905.

\end{thebibliography}
\end{document}